\documentclass{itor}

\usepackage{natbib}
\usepackage[figuresright]{rotating}

\usepackage{url}
\usepackage{algorithmic}
\usepackage{algorithm}

\usepackage{amsmath,amsthm}
\newtheorem{theorem}{Theorem}
\newtheorem{lemma}{Lemma}

\theoremstyle{definition}

\newtheorem{remark}{Remark}

\theoremstyle{remark}

\jname{International Transactions in Operational Research}
\jvol{XX}
\jyear{20XX}
\doi{xx.xxxx/itor.xxxxx}

\usepackage{amsmath,amsfonts,amsthm,amssymb,verbatim,graphicx,subfigure,color}

\usepackage{tikz}

\def\k{{\mathbf k}}
\def\u{{\mathbf u}}
\def\t{{\mathbf t}}
\newcommand{\uno}{{\mathbf{1}}}

\newcommand{\Z}{\mathbb{Z}}

\newtheorem*{theorem*}{Theorem}

\newtheorem{corollary}{Corollary}

\begin{document}

\title{A unified approach for domination and packing problems in graphs. 
}

\author[Running Author]{E. Hinrichsen\affmark{a}, G. Nasini\affmark{a,b} and N. Vansteenkiste\affmark{a,b,$\ast$}}

\affil{\affmark{a}Depto. de Matem\'atica (FCEIA), Universidad Nacional de Rosario, 
Rosario, Argentina}
\affil{\affmark{b}CONICET, Rosario, Argentina}

\email{ericah@fceia.unr.edu.ar [E. Hinrichsen]; nasini@fceia.unr.edu.ar [G. Nasini];\\ natali@fceia.unr.edu.ar [N. Vansteenkiste]}

\thanks{\affmark{$\ast$}Author to whom all correspondence should be addressed (e-mail: natali@fceia.unr.edu.ar).}

\historydate{Received DD MMMM YYYY; received in revised form DD MMMM YYYY; accepted DD MMMM YYYY}

\begin{abstract}

In this paper, we introduce new concepts of domination and packing functions in graphs, which generalize, respectively, the labelled dominating and packing functions defined by Lee and Chang in 2008, and Hinrichsen et al. in 2019.
These generalized functions offer a unified and simpler framework for addressing many of the variations of domination and packing concepts in graphs explored in the literature. Interestingly, their associated optimization problems turn out to be equivalent, 
providing insight 
to explain the observed coincidences in computational complexity results for graph classes where both problems, the domination one and its corresponding packing variation, have been analyzed.
This equivalence also allows us to solve some computational complexity open questions, for some graph classes.

Furthermore, we prove that the generalized problems remain solvable in polynomial time for graphs with bounded clique-width and strongly chordal graphs.

\end{abstract}


\maketitle

\section{Introduction}

 Domination and packing problems provide effective models for various operations research applications, particularly utility location optimization. The classical domination problem has been extensively studied, resulting in multiple variations inspired by practical applications (see \cite{Haynes} for a comprehensive survey). These variations often incorporate additional constraints tailored to specific contexts, further enriching graphs' theoretical framework of domination.
  The concept of 2-packing, a set of vertices such that no closed neighborhood contains more than one vertex from the set, serves as the \emph{packing counterpart} for the domination sets. 
 For many variations of the domination problem, their \emph{reverse}  ones, introduced in this paper as \emph{dual domination-packing pair} of problems, have been proposed in the literature and studied regarding their computational complexity in various classes of graphs.

 \cite{leechang} introduces a domination problem encompassing most of these variations as special cases. Its dual problem, studied by Hinrichsen et al., also generalizes the packing problems presented in the literature. 

In this work, we propose a new 
dual domination-packing pair of problems 
that generalizes labelling problems and simplifies their definitions. Specifically, in this generalized framework, instances are characterized by a graph and two non-negative vectors encoding bounds on the function values and on the accumulated quantities within each closed neighborhood. As a result, optimization problems previously addressed by \emph{labelling concepts} can now be viewed as particular cases of this more simplified formulation.

Interestingly, these two newly defined problems have proven to be equivalent, providing insights into the observed coincidences in the computational complexity of pairwise dual problems for specific graph classes, results
that were previously derived independently.

The structure of this paper is as follows. The next section presents the notation, definitions, and preliminary results required to understand the topics discussed. Subsequently, Section \ref{unificando} introduces the results that unify and simplify the definitions of labelling problems, establish the duality between certain variations of packing and domination problems, and clarify the relationships among them. This section also includes the definition of the new dual domination-packing pair of problems, which generalizes the labeling problems. In Section \ref{complexity}, the known computational complexity results from the literature are, presented, and the equivalence of the generalized problems is proved. Building on this equivalence, Section \ref{graphclasses} extends computational complexity results for packing problems on dually and doubly chordal graphs and demonstrates that the generalized problems remain polynomial-time solvable for graphs with bounded clique-width and strongly chordal graphs. Finally, Section \ref{conclusiones} presents the conclusions and discusses open problems.

\section{Preliminaries} \label{preliminary}

Given $n,m\in \Z_+$ with $n\leq m$ and $a\in \Z$, we denote $[n,m]=\{j\in \Z_+: n\leq j\leq m\}$, $[n]=[1,n]$ 
and  $a^+=\max\{a,0\}$. We assume that $\max \emptyset = - \infty$.  The vector $\mathbf 1 \in \Z_+^n$ (resp. $\mathbf 0\in \Z_+^n$) has all its components equal to 1 (resp. to 0).
For a finite set $V$, a function $f:V\rightarrow \Z$ can be also considered as a vector $\mathbf f \in \Z^V$. For  $U\subset V$, $f(U)$  denotes $\sum_{u\in U} f(u)$.

In this work, every graph is simple and undirected. When the sets of vertices and edges of a graph $G$ are not given, they are denoted by $V(G)$ and $E(G)$, respectively. Given $n\in \Z_+$, $K_n$, $S_n$, and $C_n$ ($n\geq 3$) are, respectively, the complete graph, the edgeless graph, and the cycle
on $n$ vertices.  

Given graphs $G$ and $H$ and a vertex $v$ of $G$, the graph obtained \emph{replacing} $v$ \emph{by} $H$ is the graph with vertices in $V(G)\setminus \{v\})\cup V(H)$ and whose edges are those in $E(G [V\setminus \{v\}])\cup E(H)$ together with those connecting every vertex of $H$ with every vertex of $N_G(v)$.

A \emph{clique} of a graph $G$ is a subset of pairwise adjacent vertices. while a \emph{stable set}, is a subset of pairwise nonadjacent vertices.
A \emph{split graph} is a graph such that its vertex set may be partitioned into a stable set and a clique.

A  \emph{chord} in a cycle is an edge connecting two non consecutive nodes in the cycle. 
A graph $G$ is \emph{chordal} if every cycle of length at least 4 has a chord.

For $v\in V(G)$, $N_G(v)$ denotes the open neighborhood of $v$ in  $G$ and $N_G[v]=N_G(v)\cup \{v\}$ is its closed neighborhood.
When the graph is clear from the context, we simply denote $N[v]$ and $N(v)$.
The \emph{degree of} $v$ is $\delta(v)= |N_G(v)|$ and $\Delta(G)=\max \{\delta(v): v\in V(G)\}$. 

A vertex $v$ is  a \emph{pendant vertex} if  $|N_G(v)|=1$ and it is called \emph {simplicial} if $N[v]$ is a clique in $G$. A vertex $u \in N[v] $ is a \emph{maximum neighbor of} $v$ if for each  $w \in N[v] $, $ N[w] \subseteq N[u]$.

Given a graph $G$ with $n$ vertices and an order $(v_1, \ldots, v_n)$ of its vertices, let $G_i=G[v_i, \ldots, v_n]$, $i\in [n]$. The order is a  \emph {strong elimination order of} $G$ if, for every $i\in [n]$, it holds that
$v_i$ 
is simplicial in $G_i$ and,   
for $j,k\in [n]$ such that $i \leq j \leq k $ and $ \{v_j, v_k\} \subseteq N_{G_i} [v_i] $ it holds that $ N_ {G_i}[v_j] \subseteq N_{G_i} [v_k]$.
Additionally, $(v_1, \ldots, v_n)$ is an \emph{order of maximum neighborhoods of }$G$ if, for $i \in [n] $, $ v_i $ has a maximum neighbor in $ G_i$.  

A graph is \emph{strongly chordal} if it admits a strong elimination ordering and it is \emph {dually chordal} if it admits an order of maximum neighborhood \cite{bran}. 
A graph is \emph {doubly chordal} if it is chordal and dually chordal. 
Several characterizations of strongly chordal graphs are presented in \cite{farber} and one of them yields a polynomial recognition algorithm and constructs strong elimination orderings.

The \emph{clique-width} of a graph 
is defined as the minimum number of labels needed to construct $G$ using the following four operations: creation of a new vertex $v$ with label $i$,  
disjoint union of two labelled graphs, 
join by an edge every vertex labelled $i$ to every vertex labelled $j$, $i\neq j$, and  
renaming label $i$ to label $j$.  
Some important classes of graphs have clique-width bounded by a constant. In particular,  cographs are the class of graphs with clique-width at most 2, trees and distance hereditary graphs have clique-width at most 3 and $P_4$-tidy graphs have clique-width at most 4. Many optimization problems that are NP-hard for more general graph classes may be solved
efficiently on graphs of clique-width bounded by a constant \citep{cmr}.

\begin{remark} \label{pendientes}
In what follows, it will be important to note that adding pendant vertices to a graph preserves the conditions of being chordal, dually chordal, or bounded clique-width.
    \end{remark}

Let us now introduce the problems we address in this paper. As mentioned before, the 
\emph{Labelled Domination problem} ($\mathcal L$dom) is presented 
 in \cite{leechang}.

A \emph{labelling function} $L$ assigns to each vertex $v$ of a graph a pair $L(v) = (t(v), k(v))$ such that 
$k(v) \in \Z$ and
$t(v)=F$ (meaning that $v$ is a free vertex) or $t(v)\in Y=\{I+t d: t\in [0,\ell]\}$,  for given $I\in \Z$ and $\ell, d \in \Z_+$.  Observe that a labelling function $L$ of a graph $G=(V,E)$ can be given by integers $I\in \Z$,  $\ell, d \in \Z_+$ and vectors $\mathbf t\in (Y \cup \{F\})^V$  
and $ \k\in \Z^V$. So, we denote it as $L=(I,d,\ell, \mathbf t,\k)$. 

Given a labelling function $L=(I,d,\ell, \mathbf t,\k)$, a function $f : V \rightarrow Y$ is an $L$-\emph{dominating function of} $G$ if, 
for all $v\in V$, 
 $f (N_G[v])\geq k(v)$ and, if $t(v) \neq F$, then $f (v) = t(v)$. The weight of $f$ is $f(V)$.  
The goal of $\mathcal L$dom is to calculate the minimum weight of an $L$-dominating function of $G$, denoted by $\gamma_L(G)$. 

The authors in \cite{leechang} prove that $\mathcal L$dom generalizes the following problems whose original definitions are presented in Appendix \ref{apendice}: Domination,  $k$-tuple domination ($k$LD), $\{k\}$-dominating functions ($\{k\}$DF), 
Signed domination, Minus domination, 
and Fault Tolerant domination. 
Other variations defined later are also generalized by $\mathcal L$dom:  M-domination  (Mdom) - also presented in \cite{LNreductions} as Generalized Multiple Domination problem (GMD) - and Multiple Domination (LD).  

For most of the domination problems mentioned above, a \emph{dual packing problem} was defined in the literature.
In particular, we can mentioned the following: 2-Packing, $k$-Limited Packing ($k$LP), $\{k\}$-packing  ($\{k\}$PF), 
Generalized Limited Packing (GLP), and  Limited Packing (LP)
(see Appendix \ref{apendice} for original definitions and references). As in the case of $\mathcal L$dom, in Remark 4 in \cite{hls} it  is shown that all the packing problems mentioned above can also be thought as particular instances of the \emph{Labelled Packing problem} 
defined as follows: 

Given a labelling function $L=(0,1,\ell, \mathbf t,\k)$ of $G=(V,E)$, a function  
$f:V\rightarrow [0,\ell]$ is an $L$-\emph{packing function of} $G$ if, for all $v\in V$, $f(N[v])\leq k(v)$ and if $t(v)\neq F$, then $f(v)=t(v)$. The objective of the Labelled Packing problem is to calculate the minimum weight of an $L$-packing function of $G$.

\section{Unifying and generalizing domination and packing problems} \label{unificando}

Let us begin by presenting some results that will allow us to more clearly establish the \emph{domination-packing duality} between the problems mentioned in the previous section.
  
Theorem 1 and Theorem 3 in \cite{ALT}  
prove that to solve $\mathcal L$dom, it is enough to consider instances with $I= 0$ and $d = 1$. 
Indeed, to obtain $\gamma_L(G)$ for a graph $G$ with $n$ vertices and $L=(I,d,\ell, \mathbf t, \k)$, we can calculate $\gamma_{L^*}(G)$ for $L^*=(0,1,\ell, \mathbf t, \k^*)$ where $k^*(v)=\lceil\frac{k(v)-I(\delta(v)+1)}{d}\rceil$ for all $v$. It follows that
$\gamma_L(G)=d \, (\gamma_{L^*}(G)+I n)$.  

 Observe that for labelling functions with $I=0$ and $d=1$, it holds that $Y=[0,\ell]$, which motivated the definition of its dual packing problem by focusing exclusively on these types of instances.
 Moreover, 
 Proposition 5 in \cite{hls} shows that the Labelled Packing Problem can be reduced to a labelling functions with $I=0$, $d=1$ and $\mathbf t\in \{0,F\}^V$.  
Thus, if  
\begin{center}
$L_{\ell, \k, \mathbf t}(G)=\min \{f(V): f \text{ is an } (0,1,\ell,\k,\mathbf t)\text{-packing  function of } G\}$    
\end{center}
for given $\ell\in \Z_+$, $\k\in \Z_+^{V}$, and $\mathbf t\in \{0,F\}^V$,  we define $\mathcal L$PF as follows: 

\medskip
{$\mathcal L$\bf{PF}: Labelling Packing Function Problem} 

Instance: A graph $G=(V,E)$, $\ell\in \Z_+$, $\k\in \Z_+^{V}$, $\mathbf t\in \{0,F\}^V$.

Objective: To obtain $L_{\ell, \k,\t}(G)$. 

\medskip 

The following lemma shows that instances of $\mathcal L$dom can be further restricted to those with $I= 0$, $d = 1$, and $\mathbf t \in \{0, F\}^V$. 
 
\begin{lemma} \label{W0}
Let be an instance of $\mathcal L$dom given by $G=(V,E)$ and $L=(0,1, \ell, \k, \mathbf t)$ with $\ell \in \Z_+$, $\mathbf t\in ([0,\ell] \cup \{F\})^V$, and $ \k\in \Z_+^V$.
Let $\mathcal F=\{v\in V: t(v)=F\}$ and $L_0=(0,1,\ell, \mathbf{k_0}, \mathbf {t_0})$ such that,
for $v\in V$,  $k_0(v)=(k(v)-t(N[v]\setminus \mathcal F))^+$, $t_0(v)=0$ if $t(v)\neq F$ and $t_0(v)=F$,  otherwise. 
Then, $\gamma_L(G)=\gamma_{L_0}(G)+t(V\setminus \mathcal F)$
\end{lemma}

\begin{proof}

Given  $f$ an $(0,1,\ell,\k,\mathbf t)$-dominating function of $G$ we define $f':V \rightarrow [0,\ell]$ such that if $w\in \mathcal F$, $f'(w)=f(w)$,  and $f'(w)=0$ otherwise.
Let us see that $f'$ is a $(0,1,\ell,\mathbf{k_0},\mathbf{t_0})$-dominating function.

For all $v\in V$ it holds that
$$f(N[v])=f( N[v]\cap \mathcal F)+
f(N[v]\setminus \mathcal F)
=f(N[v]\cap \mathcal F)+
t(N[v]\setminus \mathcal F)\geq k(v).$$ 
Then, $f(N[v] \cap \mathcal F)  \geq k(v)-t(N[v]\setminus \mathcal F)$. Since by definition of $f$, $f(N[v] \cap \mathcal F)  \geq 0$, we obtain  
$$f'(N[v]) \geq f(N[v] \cap \mathcal F) \geq  (k(v)-t(N[v]\setminus \mathcal F))^+= k'(v).$$
We have proved that $f'$ is a $L_0$-dominating function of $G$ with $f'(V)=f(V)-t(V \setminus \mathcal F)$.

Conversely, given an $(0,1,\ell,\mathbf{k_0},\mathbf{t_0})$-dominating function  $f'$ of $G$ we define $f(v)=t(v)$ if $v\notin \mathcal F$ and $f(v)=f'(v)$ otherwise. It holds that, for all $v\in V$,
$$f(N[v])
=f'(N[v]\cap \mathcal F)+
 t(N[v]\setminus \mathcal F)
\geq k'(v) + t(N[v]\setminus \mathcal F)\geq k(v).$$
Then, $f$ is a  $(0,1,\ell,\mathbf{k},\mathbf{t})$-dominating function of $G$ with $f(V)=f'(V)+t(V \setminus \mathcal F)$.
\end{proof}

\medskip 

As mentioned, the result above allows us to consider only those instances of $\mathcal L$dom with $I=0$, $d=1$, and $\mathbf t\in \{0,F\}^V$. 
We denote this problem with $\mathcal L$DF.
Formally, given $\ell\in \Z_+$, $\k\in \Z_+^{V}$, and $\mathbf t\in \{0,F\}^V$, let 
\begin{center}
    $\gamma_{\ell, \k, \mathbf t}(G)=\min \{f(V): f \text{ is an } (0,1,\ell, \k,\mathbf t)\text{-dominating function of } G\}$. 
    \end{center}
    Then, we define:

\medskip 

$\mathcal L${\bf{DF: Labelling Dominating Function Problem}}

Instance: A graph $G=(V,E)$,  $\ell \in \Z_+$, $\k\in \Z_+^{V}$, $\mathbf t\in \{0,F\}^V$.

Objective: To obtain 
$\gamma_{\ell, \k, \mathbf t}(G)$.

\medskip 

Some computational complexity results consider $\mathcal L$DF and $\mathcal L$PF restricted to instances with a fixed given $\ell \in \Z_+$. We call them $\ell$-$\mathcal L$DF and $\ell$-$\mathcal L$PF, respectively. 
We can establish the following:
\begin{remark}
 
GLP and Mdom are respectively equivalent to  $1$-$\mathcal L$PF and  $1$-$\mathcal L$DF.
    Indeed, the equivalence between GLP and $1$-$\mathcal L$PF directly follows from their definitions. 
    In the other side, $M$dom is  defined in \cite{liaochang} as $\mathcal L$dom  restricted to instances with  $I=0$, $d=\ell=1$, and $\mathbf t\in \{1, F\}^V$.  Applying Lemma \ref{W0} we obtain  that $M$dom is equivalent to 
$1$-$\mathcal L$DF. 
\end{remark}

Moreover, we also consider $\ell$-$\mathcal L$DF and $\ell$-$\mathcal L$PF restricted to instances with $t(v)=F$ for all vertex $v$, and we call them as Free $\ell$-$\mathcal L$DF and Free $\ell$-$\mathcal L$PF. With these definitions,  the following easily follows: 

\begin{remark}
The Fault Tolerant Domination problem 
corresponds to Free $1$-$\mathcal L$DF. 
Moreover, LD and LP are equivalent to Free $1$-$\mathcal L$DF and Free $1$-$\mathcal L$PF, respectively, when restricted to instances with $\k=k\uno$, for some $k\geq 2$. Similarly, $k$LD and $k$LP are equivalent to Free $1$-$\mathcal L$DF and Free $1$-$\mathcal L$PF, respectively, when restricted to instances with $\k=k\uno$, for a fixed $k\geq 1$.
Moreover, given $k\geq 1$, $\{k\}$DF and $\{k\}$PF are respectively equivalent to Free $\ell$-$\mathcal L$DF and Free $\ell$-$\mathcal L$PF when restricted to instances with $\k= k \uno$, for any $\ell \geq k$.  In addition, the Signed Domination problem is equivalent to Free $1$-$\mathcal L$DF restricted to instances with $k(v)= 1+ \lceil \frac{\delta(v)}{2} \rceil$, for all $v$. Minus Domination problem is equivalent to 
Free $2$-$\mathcal L$DF restricted to instances with $k(v)= 2+ \delta(v)$, for all $v$.
\end{remark}

We now present a new pair of dual domination-packing concepts in graphs.
$\mathcal L$DF and $\mathcal L$PF. 
Let $G=(V,E)$ be a graph, and $\k,\u$, two vectors in $\Z_+^V$.  
A function $f:V\rightarrow \Z_{+}$ is a $(\k,\u)$-\emph{dominating (resp. packing) 
function}  
of $G$ if it satisfies
$f(v)\leq u(v)$  and $f(N[v])\geq k(v)$ (resp. $f(N[v])\leq k(v)$), for all $v\in V$.
Denoting by  $\mathcal D_{\k,\u}(G)$ (resp. $\mathcal L_{\k,\u}(G)$ ) the set of all $(\k,\u)$-dominating (resp. packing) functions of $G$, the associated optimization problems  are formulated as follows:

\medskip 

{\bf{GDF}: Generalized Dominating Function Problem} 

Instance: A graph $G$, \, $\k,\u\in \Z_+^{V}$.

Objective: To obtain $\gamma_{\k,\u}(G)=\min \{f(V): f \in  \mathcal D_{\k,\u}(G)\}$.

\medskip 

{\bf{GPF}: Generalized Packing Function Problem} 

Instance: A graph $G$,  \, $\k,\u\in \Z_+^{V}$.

Objective: To obtain $L_{\k,\u}(G)=\max \{f(V): f \in  \mathcal L_{\k,\u}(G)\}$.

\begin{remark}   \label{generaliza}
Observe that, according to the definitions, $\mathcal L$DF and $\mathcal L$PF are respectively equivalent to GDF and GPF restricted to instances with $\u\in \{0,\ell \}^V$, for $\ell \in \Z_+$. Then, GDF and GPF generalize all domination and packing problems presented in the previous section.
\end{remark}

\begin{figure}[h]
\centering
\includegraphics[scale=0.5]{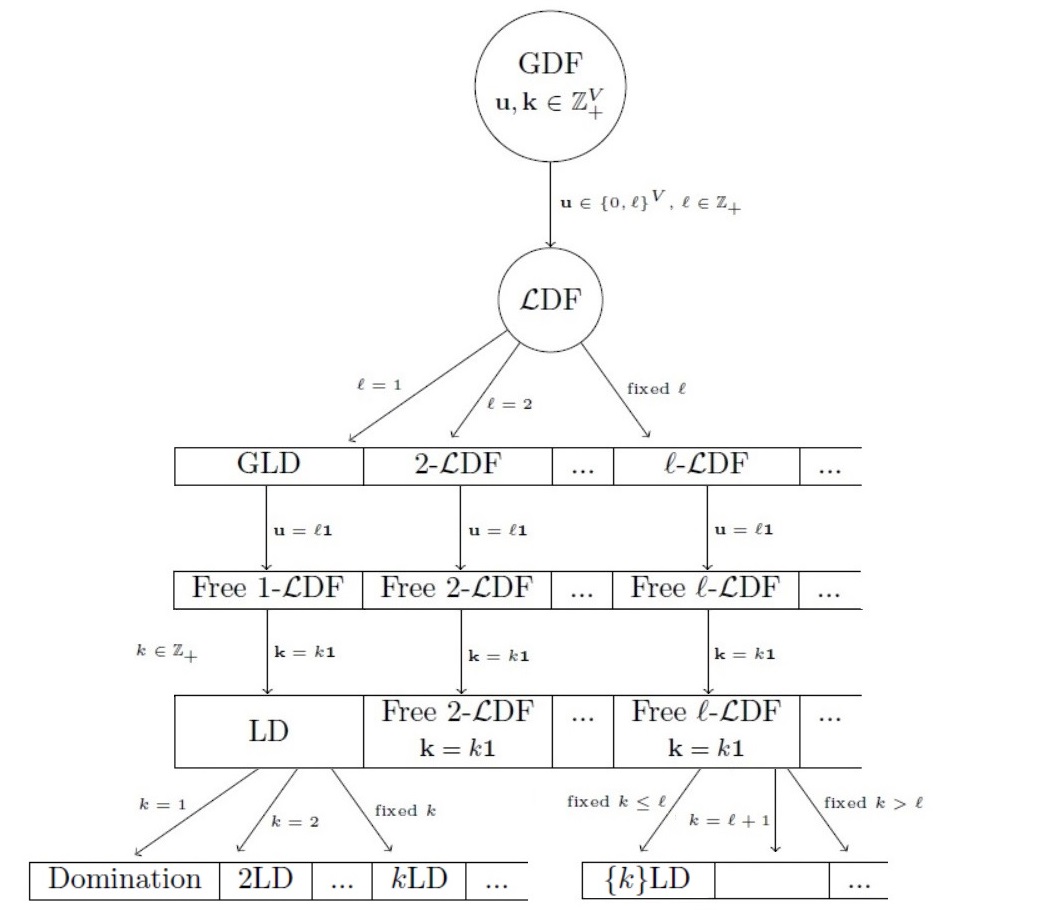}
\caption{Relationships among domination problems. The problem at the end of an arrow corresponds to particular instances of the one at its origin. The restrictions defining such special instances are presented near the arrow. }
    \label{fig: graficoambos}
\end{figure}

Given a graph $G=(V,E)$, the relationships among the domination problems analized in this work are represented in Figure \ref{fig: graficoambos}, where the acronym of each one is written in the box. Replacing the letter D with P and changing Domination to 2-packing in the figure, we obtain the relationship among packing problems. 

\medskip 
Let us make some considerations about particular conditions that can be assumed for GDF and GPF instances,  denoted by $(G,\k,\u)$. 

\begin{remark} \label{kuD}
Regarding GDF, it is clear that if $k(v)> u(N[v])$ for some $v\in V$, the instance $(G,\k,\u)$ is not feasible. 
In addition, given $f \in \mathcal D_{\k,\u}(G)$ such that $f(v)\geq  \max \{k(w): w \in N[v]\}+1$  for some $v\in V$, defining  $\tilde f(v)= f(v)-1$ and $\tilde f(w)=f(w)$ in other case,  $\tilde f \in \mathcal D_{\k,\u}(G)$  and $\tilde f(V)\leq f(V)$.  Therefore if $f \in \mathcal D_{\k,\u}(G)$ with $f(V)=\gamma_{\k,\u}(G)$,  $f(v)\leq  \max \{k(w): w \in N[v]\}$. Then, an instance $(G,\k,\u)$ of GDF is equivalent to instance $(G, \tilde\k,\tilde\u)$ where $\tilde k(v)=\min \{k(v), u(N[v])\}$,  $\tilde u(v)=\min \{u(v), \max \{k(w): w\in N[v]\}\}$ for all $v\in V$. We can then assume without loss of generality that the instances $(G,\k,\u)$ of GDF verify $k(v)\leq u(N[v])$ and  $u(v)\leq  \max \{k(w): w\in N[v]\}$, for all $v\in V$.
\end{remark}

\begin{remark} \label{kuP}
With regards to GPF, it is clear that every $f \in \mathcal L_{\k,\u}(G)$ satisfies $f(v)\leq f(N[v]) \leq k(v)$ for all $v\in V$. 
Moreover,
$f(N[v])\leq u(N[v])$ for all $v\in V$. 
Then, an instance $(G,\k,\u)$ of GPF is equivalent to instance $(G, \tilde\k,\tilde\u)$ where $\tilde k(v)=\min \{k(v), u(N[v])\}$,  $\tilde u(v)=\min \{u(v), k(v)\}$ for all $v\in V$. We can then assume without loss of generality that the instances  $(G,\k,\u)$ of GPF verify  $u(v)\leq  k(v)\leq  u(N[v])$, for all $v\in V$.
\end{remark}

\section{General Computational Complexity results} \label{complexity}

In Figure \ref{tabla} we present the best previously known result in the literature on computational complexity for the problems considered in this paper and the consequences derived from the relationships among the optimization problems and the graph classes.  
References for these results are presented in the Appendix. 
 
\begin{figure} [h]
\centering
\begin{tabular}{ l  c  c  c  c  c  c  }

   &\footnotesize{Domination}&\footnotesize $k$\bf{LD}  &  \footnotesize\bf{LD}  &  
   \footnotesize $\{k\}$\bf{DF} & \footnotesize $\ell$-\bf{$\mathcal L$DF}  & \footnotesize \bf{$\mathcal L$DF} 
   \\ 
    &\footnotesize ($2$-Packing) & \footnotesize($k$\bf{LP}) & \footnotesize(\bf{LP}) &   \footnotesize($\{k\}$\bf{PF}) & \footnotesize($\ell$-\bf{$\mathcal L$PF}) & \footnotesize(\bf{$\mathcal L$PF}) 
  \\ 
    & & \scriptsize {$k\geq 2$} &  &   \scriptsize {$k\geq 2$} &  &  
     \\
    \hline
    
   \footnotesize  bipartite  & \scriptsize NP-h  &  \scriptsize NP-h   & \scriptsize NP-h 
   &   \scriptsize NP-h   &     \scriptsize NP-h     & \scriptsize NP-h
   \\ 
     & \scriptsize (NP-h) & \scriptsize (?) & \scriptsize (NP-h) 
     
    & \scriptsize (NP-h)  & \scriptsize (NP-h) &  \scriptsize (NP-h)   
    \\ \hline
    
  \footnotesize split  & \scriptsize NP-h  &   \scriptsize NP-h &   \scriptsize NP-h    &     \scriptsize NP-h     &  \scriptsize NP-h        &  \scriptsize NP-h
  \\ 
   
 & \scriptsize (NP-h) & \scriptsize (NP-h) & \scriptsize (NP-h)  &   \scriptsize (NP-h)   & \scriptsize (NP-h)   & \scriptsize (NP-h) 
 \\ \hline 

   \footnotesize dually  &\scriptsize P & \scriptsize NP-h & \scriptsize NP-h   &  \scriptsize P  & \scriptsize NP-h   & \scriptsize NP-h
   \\ 
  \footnotesize chordal & \scriptsize (P)  & \scriptsize (?) & \scriptsize (?) &  \scriptsize (P) & \scriptsize (?) & \scriptsize (?)
  \\ \hline

  \footnotesize doubly  &\scriptsize P & \scriptsize NP-h   & \scriptsize NP-h &   \scriptsize P & \scriptsize NP-h & \scriptsize NP-h
  \\ 
  \footnotesize chordal &\scriptsize (P) & \scriptsize (?) & \scriptsize (?)  &  \scriptsize (P) & \scriptsize (?) & \scriptsize (?)
  \\ \hline
   
  \footnotesize strongly   &\scriptsize P & \scriptsize P & \scriptsize P &    \scriptsize P    &   \scriptsize P   & \scriptsize P  
  \\
  \footnotesize chordal  & \scriptsize (P) & \scriptsize (P)      &  \scriptsize (P) & \scriptsize (P) & \scriptsize (P) & \scriptsize (P)  
  \\ \hline
   
  \footnotesize bounded & \scriptsize P & \scriptsize  P & \scriptsize P & \scriptsize P & \scriptsize  P  & \scriptsize ?
  \\ 
   \footnotesize clique-width &\scriptsize (P) & \scriptsize (P)  & \scriptsize (?)  & \scriptsize (P) & \scriptsize (?) &\scriptsize (?) 
   \\ \hline

\end{tabular}
\caption{
Each entry in this table corresponds to a graph class and a domination.-packing pair of problems.  NP-h means that a problem is NP-hard, P indicates that it is polynomial-time solvable, and the question mark indicates that its computational complexity is unknown.} 
    \label{tabla}
\end{figure}

From the table, it can be observed that there is no known graph family where one of the problems of a dual domination-packing pair of problems is polynomial-time solvable and the other, NP-hard. 
The following result sheds some light on this observation:

\begin{theorem}
    \label{GPFigualGDF}
Let $(G, \k, \u)$ be an instance of GDF or GPF, with $G=(V,E)$. Let $\k'$ such that $k'(v) =u(N[v])- k(v)$ for each $v\in V$. Then, $L_{\k,\u}(G) = u(V)-\gamma_{\k',\u}(G)$ and $\gamma_{\k,\u}(G) = u(V)-L_{\k',\u}(G)$.
\end{theorem}
\begin{proof}
By Remarks \ref{kuD} and \ref{kuP} we can consider that $k(v) \leq u(N[v])$ for all $v\in V$. Then, $\k'\geq \mathbf 0$.

Given $f \in \mathcal L_{\k,\u}(G)$ such that $f(V)=L_{\k,\u}(G)$, we define $g:V \rightarrow \Z_{+}$ such that $g(v)=u(v)-f(v)$. Since $0\leq f(v)\leq u(v)$ it holds that $0\leq g(v)\leq u(v)$. Moreover, $$g(N[v])=u(N[v])-f(N[v])\geq u(N[v])-k(v) = k'(v).$$ 
Therefore, $g \in \mathcal D_{\k',\u}(G)$ with $g(V) = u(V )-f(V) = u(V)-L_{\k,\u}(G)$, implying that  $\gamma_{\k',\u}(G) \leq u(V)-L_{\k,\u}(G)$ or, equivalently,  $L_{\k,\u}(G)\leq u(V) - \gamma_{\k',\u}(G)$.

Conversely, given $g\in \mathcal D_{\k',\u}(G)$ such that $g(V)=\gamma_{\k,\u}(G)$, we define $f:V \rightarrow \Z_{+}$ with $f(v) = u(v)-g(v)$. 
It holds that $0 \leq f(v) \leq u(v)$ and
$$f(N[v]) = u(N[v])-g(N[v])\leq u(N[v])-k'(v) = k(v).$$ 
Then, $f \in \mathcal L_{\k,\u}(G)$ 
and $L_{\k,\u}(G) \geq f(V)= u(V)-\gamma_{\k',\u}(G)$. Finally, 
$$L_{\k,\u}(G) = u(V)-\gamma_{\k',\u}(G).$$

Applying the previous result to the instance $(G,\k',\u)$ of GPF, we obtain $L_{\k',\u}(G) = u(V)-\gamma_{\k,\u}(G)$ and then $\gamma_{\k,\u}(G) = u(V)-L_{\k',\u}(G)$.
\end{proof}

Theorem \ref{GPFigualGDF} it allows to easily prove the following result: 
\begin{theorem} \label{iguales}
$\mathcal L$DF and $\mathcal L$PF, as well as   $\ell-\mathcal L$DF and $\ell-\mathcal L$PF,  Free $\ell-\mathcal L$DF and Free $\ell-\mathcal L$PF, 
have the same computational complexity for any family of graphs. 
\end{theorem}
Regardless of Theorem 1, the equivalence between 1-LDF and 1-LPF for any family of graphs can easily be obtained by applying Lemma 1 before Proposition 1 in \cite{LNreductions}.

Observe that no such direct consequence can be stated on the dual pairs $k$LD--$k$LP and $\{k\}$DF--$\{k\}$PF. Indeed, on these problems we have $\k= k \uno$ and $\u= \ell \uno$, for $\ell=1$ and $\ell=k$, respectively. Then, the values of vector $\k'$ in Theorem \ref{GPFigualGDF} become  $k'(v)=u(N[v])-k(v)=\ell (\delta(v)+1)- k$ , depending on the degree of each vertex $v\in V$. 
Therefore,
we can establish:

\begin{corollary} \label{cor2}
$k$LD  and $k$LP, as well as $\{k\}$DF and $\{k\}$PF,  have the same computational complexity for any family of regular graphs. 
\end{corollary}

No specific results on computational complexity for Free $\ell$-$\mathcal L$DF and Free $\ell$-$\mathcal L$PF are known. The following lemma proposes a polynomial-time reduction that allows us to obtain some results in this direction: 

\begin{lemma} \label{reduccionafree}
Given $\ell \in \Z_+$, $\ell$-$\mathcal L$PF can be reduced to Free $\ell$-$\mathcal L$PF in polynomial time.
\end{lemma}
\begin{proof}
Let $(G,\k,\u)$ be an instance of $\ell$-$\mathcal L$PF with $G=(V,E)$. Let $\tilde G$  be the graph obtained by adding a pendant vertex $v'$ to each $v\in V$ such that $u(v)=0$. Consider $(\tilde G, \tilde \k, \uno)$, the instance of Free $\ell$-$\mathcal L$PF  where $\tilde k(v)=k(v)$ for all $v\in V$, and $\tilde k(v)=0$, otherwise. 
\\ Since $\tilde k(w)=0$ for all $w\in V(\tilde G)\setminus V$, every $\tilde f\in \mathcal L_{\tilde \k, \tilde \u}(\tilde G)$ verifies that $\tilde f(v)=0$ for all $v\in V$ such that $u(v)=0$. Then, if $f$ is defined with domain $V$ such that $f(v)=\tilde f(v)$ for all $v\in V$, it holds that $f\in \mathcal L_{\k,\u}(G)$.

Conversely, for every $f\in \mathcal L_{\k,\u}(G)$ the function define as $\tilde f(v)=f(v)$ for all $v\in V$ and $\tilde f(v)=0$ for all $v\in V(\tilde G)$ verifies that $\tilde f\in \mathcal L_{\tilde \k,\tilde \u}(\tilde G)$. 
Clearly, $\tilde f(V(\tilde G))=f(V)$. 
\end{proof}

From the previous lemma and Theorem \ref{iguales} we obtain:

\begin{corollary} \label{pendientefree} 
     Let $\mathcal F$ be a family of graphs closed under addition of pendant vertices and $\ell\in \Z_+$. If $\ell$-$\mathcal L$PF (resp. $\ell$-$\mathcal L$DF ) is NP-hard for $\mathcal F$ then 
      Free $\ell$-$\mathcal L$PF (resp.  Free $\ell$-$\mathcal L$DF) is NP-hard for $\mathcal F$. In addition, if Free $\ell$-$\mathcal L$PF (resp.  Free $\ell$-$\mathcal L$DF) is polynomial-time solvable for $\mathcal F$ then $\ell$-$\mathcal L$PF (resp. $\ell$-$\mathcal L$DF ) is polynomial-time solvable for $\mathcal F$.
\end{corollary}

\section{Computational Complexity results} \label{graphclasses}

In this section, we complete some entries in Table \ref{tabla} and analyze the computational complexity of the new problems for the graph classes considered there. From Remark \ref{generaliza} we easily obtain that GDF and GPF are NP-hard in bipartite and split graphs. 
\subsection{On doubly and dually chordal graphs}
It can be observed that doubly and dually chordal graphs are the only graph classes considered in Table \ref{tabla} where one of the problems is polynomial-time solvable ($\{k\}$DF) and all the others are NP-hard. Regarding their packing dual counterparts, only the computational complexity of $\{k\}$PF  is currently known.
From Remark \ref{generaliza}
 and Theorem \ref{iguales} we can already establish:
\begin{theorem}
For a given $\ell \in \Z_+$,  $\ell$-$\mathcal L$PF, and  $\mathcal L$PF are NP-hard on dually and doubly chordal graphs. Therefore, GDF and GPF are also NP-hard for these graph classes. 
\end{theorem}

The computational complexity of $k$LP or LP on doubly chordal graphs can not be obtained from Theorem \ref{iguales}. However, we can prove:

\begin{theorem}
LP is NP-hard on doubly chordal graphs.    
\end{theorem}

\begin{proof}
Let  $(G,\k,\uno)$ be an instance of  Free $1$-$\mathcal L$PF with $G=(V,E)$ and 
$k^*=\max \{k(v): v\in V\}$. Observe that 
by Remark \ref{kuP}
we can consider $k^* \leq  \Delta(G)+1$. 
Consider the graph $G'=(V',E')$ obtained by adding a set $S_v$ of $k^*-k(v)$ pendant vertices to every $v\in V$, and let $\k'=k^*\uno\in \Z^{V'}$. Clearly,  
$(G',\k',\uno)$ is an instance of LP.

Let $f\in \mathcal L_{\k,\uno}(G)$. It is not difficult to see that $f':V' \rightarrow \{0,1\}^{V'}$ such that $f'(v)= f(v)$ if $v\in V$ and $f'(v)=1$, otherwise, satisfies $f'\in \mathcal L_{\k',\uno}(G')$. Moreover, since $f'(V')=f(V)+\sum_{v\in V}(k-k(v))$ we obtain $L_{\k',\uno}(G')\geq L_{\k,\uno}(G)+\sum_{v\in V}(k^*-k(v))$.

For the converse, let us first observe that for any function $g\in \mathcal L_{\k',\uno}(G')$ we can obtain a function $\tilde g\in \mathcal L_{\k',\uno}(G')$ such that $\tilde g (V')\geq g(V')$ and $\tilde g(w)=1$ for any pendant vertex $w\in V'$. 
Indeed, let $w$ be a pendant vertex of $G'$  such that $g(w)=0$ and $v\in N(w)$. If $g(z)= 1$ for some $z\in N[v]\setminus \{w\}$, we can define $\tilde g$ interchanging the values of $w$ and $z$. Otherwise, since $\k'\geq \uno$, we can define $\tilde g(w)=1$ without changing any other value of $g$.

Then, let $f'\in \mathcal L_{\k',\uno}(G')$ such that $f'(V')=L_{\k',\uno}(G')$ and $f'(w)=1$ for  all $w\in S_v$ and $v\in V$. 
Let $f$ be the projection of $f'$ on $V$. Clearly, 
$$f(N_{G}[v])=f'(N_{G'}[v])-(k^*-k(v))\leq k^*-(k^*-k(v))=k(v).$$
Therefore, $f\in \mathcal L_{\k,\uno}(G)$ and  
$$ L_{\k,\uno}(G) \geq f(V)=
L_{\k',\uno}(G')-\sum_{v\in V}(k-k(v)),$$  
proving that $L_{\k',\uno}(G')= L_{\k,\uno}(G)+\sum_{v\in V}(k^*-k(v))$. 

The transformation is linear-time and closed on dually and doubly chordal graphs (Remark \ref{pendientes}).
\end{proof}

As well as in bipartite graphs, the computational complexity of $k$LP on dually and doubly chordal graphs is unknown.

\subsection{In strongly chordal graphs} \label{schordal}
Despite $k$LD and LP are NP-hard on doubly chordal graphs,  GDF and GPF remain polynomial-time solvable on the subclass of strongly chordal graphs.
 
\begin{theorem}
Given a strong elimination ordering of vertices of $G=(V,E)$ GDF and GPF on any instance $(G,\k,\u)$ can be solved in $O(|V|+|E|)$-time. 
\end{theorem}

\begin{proof}
We prove that Algorithm \ref{alg} solves GPF for any instance $(G, \k,\u)$.
 
\begin{algorithm} [h] 
\begin{algorithmic} [1]

\STATE{\emph{Input:}}
$G=(V,E)$, $\k,\u\in \Z_+^V$, $(v_1,v_2,\ldots,v_n)$ strong elimination ordering of $G$.

\STATE{\emph{Initialization:}} for $i=1$ to $n$ do $f(v_i):=0$.

\STATE{\emph{Augmenting the weight of the function:}}
\hspace{.2cm} for $i=1$ to $n$ do 
\STATE
\hspace{.5cm} $M:=\min\{k(v)- f(N[v]):v\in N[v_i] \}$ 
\STATE
\hspace{.5cm}  
 $f(v_i):=min\{ M, u(v_i)\}$.

\STATE{Return function $f$}
   
\end{algorithmic}
\caption{}
\label{alg}
\end{algorithm}

For $i\in [0,n]$, let  $f^{(i)}$ be the function $f$ obtained at the end of the $i$-th step of the algorithm. Therefore  $f^{(0)}$ is the null function and $f^{(n)}$ is the output of the algorithm. Let $f=f^{(n)}$.

Observe that, for $i\in [0,n-1]$, $f^{(i+1)}$ is obtained from $f^{(i)}$ just by modifying the value on $v_i$ in step 5. 
Therefore, given $j\in [n]$, $f^{(i)}(v_j)=f^{(0)}(v_j)=0$ for $i\leq j-1$ and $f^{(i)}(v_j)= f^{(j)}(v_j)$ for $i\geq j$.  
Clearly $f^{(i)}(v)\leq u(v)$ for all $v\in V$ and $i\in[0,n]$.

We first prove by induction that  
$f^{(i)}\in \mathcal L_{\k,\u}(G)$ for all $i\in [0, n]$.  
Clearly, $f^{(0)}\in \mathcal L_{\k,\u}(G)$ and assume that  $f^{(i-1)}\in \mathcal L_{\k,\u}(G)$ for some $i \in [n]$. To prove that $f^{(i)}\in \mathcal L_{\k,\u}(G)$, it only remains to check that 
$f^{(i)}(N[v])\leq k(v)$ for all $v\in V$.

If $v\notin N[v_{i}]$, from previous observations, 
$f^{(i)}(N[v])=f^{(i-1)}(N[v])\leq k(v)$. 
In addition, if $v\in N[v_i]$, then 
$$f^{(i)}(N[v])=f^{(i-1)}(N[v]) + f^{(i)}(v_i)\leq f^{(i-1)}(N[v]) + M \leq$$$$\leq    f^{(i-1)}(N[v]) + k(v) - f^{(i-1)}(N[v])= k(v).$$ 

Then, $f^{(i)}\in \mathcal L_{\k,\u}(G)$ for all $i\in [0,n]$ and it remains to prove that $f(V)= L_{\k,\u}(G)$.

Let $h\in \mathcal L_{\k,\u}(G)$ with $h(V)= L_{\k,\u}(G)$ and such that the cardinality of  $W_h=\{v\in V: f(v)\neq h(v)\}$ is  minimum. 
We want to prove that $W_h$ is empty. 
 
Suppose that  $W_h\neq\emptyset$ and let $l$ be the minimum index such that $v_l\in W_h$, i.e. $h(v_j)=f(v_j)$ if $j<l$ and $h(v_l)\neq f(v_l)$. We analyze the two possible cases:  $h(v_l)>f(v_l)$ and $h(v_l)<f(v_l)$.

\smallskip

\noindent \textbf{Case 1} $h(v_l)>f(v_l)$. 

Since $f(v_l)=f^{(l)}(v_l)$ and $h\in \mathcal L_{\k,\u}(G)$, it holds that $f^{(l)}(v_l) < h(v_l)\leq u(v_l)$. Then $f^{(l)}(v_l)=M=\min \{k(v)-f^{(l-1)}(N[v]): v\in N[v_l])$.
Let $v_c\in N[v_l]$ such that  $M=k(v_c)- f^{(l-1)}(N[v_c])$. Then,
$$h(N[v_c])\geq h(N[v_c]\cap \{v_j: j\leq l-1\})+h(v_l)= f(N[v_c]\cap \{v_j: j\leq l-1\})+h(v_l)>$$
$$> f(N[v_c]\cap \{v_j: j\leq l-1\})+f(v_l)=f^{(l)}(N[v_c]\cap \{v_j: j\leq l-1\})+f^{(l)}(v_l) \geq $$ 
$$ \geq f^{(l)}(N[v_c]\setminus\{v_l\})+f^{(l)}(v_l).$$
Since $f^{(l)}(N[v_c] \setminus \{v_l\})=f^{(l-1)}(N[v_c])=k(v_c)-M$ and $f^{(l)}(v_l)=M$, we obtain  
$$h(N[v_c]) > f^{(l)}(N[v_c]\setminus\{v_l\})+f^{(l)}(v_l) = (k(v_c)-M)+M=k(v_c),$$ contradicting the fact that  $h\in \mathcal L_{\k,\u}(G)$.

\medskip 

\noindent \textbf{Case 2} 
 $h(v_l)<f(v_l)$.

\medskip 

Let us first prove that  it is not possible that  $h(N[v])+f(v_l)-h(v_l)\leq k(v)$ for all $v\in N[v_l]$. Indeed, if this were so, defining  $g(v_l)=h(v_l) +f(v_l) -h(v_l)=f(v_l) $ and $g(v)=h(v)$ in other case, clearly holds that  $g(v)\leq u(v)$ for all $v\in V$. Moreover, for $v\notin N[v_l]$, $g$ satisfies $g(N[v])=h(N[v])\leq k(v)$. 
Finally, if  $v\in N[v_l]$, we obtain:
$$g(N[v])=h(N[v]) - h(v_l)+g(v_l)=h(N[v])-h(v_l)+f(v_l) \leq k(v).$$
Then, $g\in \mathcal L_{\k,\u}(G)$ and $g(V)> h(V)=L_{\k,\u}(G)$, a contradiction.  
\medskip 
  
Let  $P=\{v\in N[v_l]:h(N[v])+f(v_l)-h(v_l)>k(v)\} \neq \emptyset$ and  $s$  be the minimum index among all the vertices in $P$.

Observe that, for each $v\in P$, it holds that
 $$h(N[v])+f(v_l)-h(v_l)>k(v)\geq f(N[v]).$$ Then, $h(N[v]) -h(v_l)> f(N[v]) - f(v_l)$ or, equivalently, 
 $h(N[v] \setminus \{v_l\})> f(N[v] \setminus \{v_l\})$. 
 
 This implies that, for each $v\in P$, there exists $v_i\in N[v]$  with $i>l$ such that $h(v_i)>f(v_i)$ and 
 $N[v]\cap\{v_x:  h(v_x)>f(v_x)\}\neq \emptyset$. In particular for $v=v_s$ which allows us to define $b$ as the minimum index among all the vertices in 
 $N[v_s]\cap\{v_x: \; h(v_x)>f(v_x)\}$. 
 
 \medskip

 Observe that $P\subseteq N[v_l]$  and it can be verified that  $P\subseteq N[v_b]$ by considering the following two cases: 
 
 \begin{enumerate}
 \item[-] Case (a): $s\leq l$. 
 
 In this case, $s\leq l<b$. Then, $v_l$ and $v_b$ belong to  $N_{G_s}[v_s]$ and by definition of strong elimination ordering,  $N_{G_s}[v_l]\subseteq N_{G_s}[v_b]$. 
 Since $P\subseteq N_{G_s}[v_l]$ we obtain $P\subseteq N[v_b]$.

 \item[-]  Case (b): $l<s$.  
 
 Let us consider $v_k\in P$. Since $l<s\leq k$, $v_s$ and $v_k$ belong to $N_{G_l}[v_l]$. By definition of strong elimination ordering, 
 $N_{G_l}[v_s]\subseteq N_{G_l}[v_k]$, for each $v_k\in P$. Since $v_b\in N_{G_l}[v_s]$, then $v_b\in N_{G_l}[v_k]$. In other words,  $P\subseteq N[v_b]$.
 \end{enumerate}

Let  $c=\min\{f(v_l)-h(v_l),h(v_b)-f(v_b)\}>0$.
We define $h'$ such that  
$h'(v_l)=h(v_l)+c$, $h'(v_b)=h(v_b)-c$, and $h'(v)=h(v)$ otherwise. 

\medskip
Note that since $f\in \mathcal L_{\k,\u}(G)$ and $h\in \mathcal L_{\k,\u}(G)$  it holds that $$h'(v_l)=h(v_l)+c\leq f(v_l)\leq u(v_l) \text{\; and \; } h'(v_b)=h(v_b)-c\leq u(v_b).$$
To prove that  $h'\in \mathcal L_{\k,\u}(G)$, it only remains to prove  that 
$h'(N[v])\leq k(v)$ for each $v\in N[v_l]$.

\medskip

If $v\in N[v_l]\setminus N[v_b]$, $$h'(N[v])=h'(N[v]\setminus\{v_l\})+h'(v_l)=h(N[v]\setminus\{v_l\})+h(v_l)+c=$$$$=h(N[v])+c\leq h(N[v])+f(v_l)-h(v_l).$$ Since $v\notin P$, $h(N[v])+f(v_l) -h(v_l)\leq k(v)$.

Otherwise, if $v\in N[v_l]\cap N[v_b]$, then  $$h'(N[v])=h'(N[v]\setminus\{v_l,v_b\})+h'(v_l)+h'(v_b)=$$$$=h(N[v]\setminus\{v_l,v_b\})+h(v_l)+c+h(v_b)-c=h(N[v]) \leq k(v).$$
We conclude that  $h'\in \mathcal L_{\k,\u}(G)$. Observe that, if $c=f(v_l)-h(v_l)$, then  $h'(v_l)=h(v_l)+c=f(v_l)$. Otherwise, if $c=h(v_b)-f(v_b)$, then 
$h'(v_b)=h(v_b)-c=f(v_b)$. Then, $h'(V) = h(V)= L_{\k,\u}(G)$. Moreover, 
$|\{v\in V:f(v)=h'(v)\}|\geq|\{v\in V:f(v)=h(v)\}|+1$, which contradicts the fact that $|W_h|$ is minimum.

Then, Cases 1 and 2 are not possible, which implies that $f(v_l)=h(v_l)$, a new contradiction with the definition of $l$. Therefore, $W_h=\emptyset$, $f=h$ and $f(V)= L_{\k,\u}(G)$.

\medskip 

The running time of Algorithm 1 is $O(\sum_{v\in V}|N[v]|)=O(|V|+|E|)$ and, by Theorem  \ref{GPFigualGDF}, there is also a $O(|V|+|E|)$-time algorithm for solving  GDF on strongly chordal graphs. 
\end{proof}

\subsection{In bounded clique-width graphs} \label{bcw}

In this section, we analyze the behavior of GDF and GPF on bounded clique-width graphs.   
Let us first prove the following. 

\begin{lemma} \label{reduction2}
Given $M\in \Z_+$, GDF and GPF, restricted to instances with $\u\leq M \uno$, can be reduced to $1$-$\mathcal L$DF in linear time. 

\end{lemma}
\begin{proof}
Let $(G, \k,\u)$ be an instance of GDF, with $G=(V,E)$ and $\u\leq \uno M$. 

Let $\mathcal N=\{v\in V: u(v)=0\}$ and $G'$, the graph obtained by replacing each vertex $v\in V\setminus \mathcal N$  with a complete graph $K^v$ of $u(v)$ vertices. In 
addition, let 
$\k' \in \Z_+^{V(G')}$ 
such that $k'(w)=k(w)$ if 
$w\in \mathcal N$ and 
$k'(\tilde v)=k(v)$ if 
$\tilde v \in V(K^v)$.   
Finally, we define $\tilde\u\in \{0,1\}^{V(G')}$, such that
$u'(v)=0$ for all $v\in \mathcal N$ and $u'(v)=1$, otherwise.
Let us prove that solving GDF on $(G,\k,\u)$ is equivalent to solving $1$-$\mathcal L$DF on $(G',\k',\u')$. 

Let $f\in \mathcal D_{\k,\u}(G)$. 
We define $f': V(G')\rightarrow\{0,1\}$  such that $ f'(V(K^v)) = f(v)$ for each $v\in V\setminus \mathcal N$,  and $f'(v)=0$, otherwise.

For all $v\in V(G')$ it holds that $$N_{G'}[w]=\cup_{z\in N_G[v]\setminus \mathcal N} K^z \cup (N_G[v]\cap \mathcal N).$$

Then, 
$$f'(N_{G'}[w])=\sum_{z\in N_G[v]\setminus \mathcal N} f'(K^z)+ f'(N_G[v]\cap \mathcal N)= $$
$$=\sum_{z\in N_G[v]\setminus \mathcal N} f(z)+ 0= f(N[v])\geq k(v)=k'(w).$$

We obtain $ f'\in \mathcal{D}_{\k',\u'}(G')$ and $f'(V(G'))=f(V)$.

For the converse, let  $ f'\in \mathcal{L}_{\k',\u'}(G')$. 
Now, we define $f: V\rightarrow \Z_+$ such that $f(v)=  f'(K^v)$ if $v\notin  \mathcal N$ and $f(v)=0$, otherwise. Let see that $f\in \mathcal L_{\k,\u}(G)$.

If $v\in V\setminus \mathcal N$,  $f(v)= f' (K^v)\leq |K^v| = u(v)$. Otherwise, $v\in \mathcal N$ and $f(v)=0=u(v)$. In addition, for all $v\in V$,
$$f(N_G[v])=\sum_{w\in N_G[v]\setminus \mathcal N} f(w) +0=\sum_{w\in N_G[v]\setminus \mathcal N} f'(K^w)=f'(N_{G'[v]})\geq k'(v)=k(v).$$
Then, $f\in \mathcal D_{\k,\u}(G)$ and 
$f(V)=f'(V(G'))$.

Since $|V(G')|\leq |V| \max\{u(v):v\in V\} \leq  |V| M$, the reduction can be done in linear-time.
\end{proof}

Considering that the reduction in Lemma \ref{reduction2} preserves the bounded clique-width property in graphs (Remark \ref{pendientes}) and that $1$-$\mathcal L$DF is polynomial-time solvable for this class of graphs (Corollary  \ref{iguales}),  we obtain the following:
\begin{theorem}
    Given $M\in \Z_+$, GDF and GPF restricted to instances with $\u\leq M \uno$ are polynomial-time solvable on bounded clique-width graphs. As a consequence, for a given $\ell \in \Z_+$, Free $\ell$-$\mathcal L$DF, Free $\ell$-$\mathcal L$PF and $\ell$-$\mathcal L$DF also are polynomial-time solvable on bounded clique-width graphs. 
\end{theorem}

The computational complexity of dual pairs GPF-GDF, and $\mathcal L$DF-$\mathcal L$PF for bounded clique-width graphs, is still open.

\section{Conclusions and open questions}
\label{conclusiones}

In this paper, we propose a generalized dual domination-packing pair of problems, and we prove that they are equivalent to any graph, offering a partial explanation for these coincidences in computational complexity results for dual pairs. Furthermore, these new problems provide a unified and simplified framework highlighting the \emph{domination-packing duality} inherent in several pairs of problems previously defined in the literature. This framework also facilitates the completion of findings for graph classes where only one issue in the dual domination-packing pair has been studied.
Moreover, our results extend known polynomial-time solvability results for strongly chordal and bounded clique-width graphs to these generalized problems.
The following table summarizes the complexity results obtained in this paper.

\begin{figure}[h]
\centering
\begin{tabular}{ l  c  c  c  c  c  c  }

   &\footnotesize $k$\bf{LD}  &  \footnotesize\bf{LD}  & \footnotesize $\ell$-\bf{$\mathcal L$DF}   &\footnotesize \bf{$\mathcal L$DF} &  \footnotesize\bf{GDF} &\footnotesize\bf{GDF} 
   \\ 
    & \footnotesize($k$\bf{LP}) & \footnotesize(\bf{LP}) &    \footnotesize($\ell$-\bf{$\mathcal L$PF}) &  \footnotesize(\bf{$\mathcal L$PF}) &  \footnotesize\bf{(GPF)} &  \footnotesize\bf{(GPF)}
  \\ 
     & \scriptsize {$k\geq 2$} &  &   &   & \scriptsize {$\u\leq M\uno$} &  
     \\
    \hline

   \footnotesize  bipartite  &\scriptsize NP-h & \scriptsize NP-h 
   &   \scriptsize NP-h   & \scriptsize NP-h & \scriptsize NP-h $^{**}$ & \scriptsize NP-h $^{**}$
   \\ 
    & \scriptsize (?)  & \scriptsize (NP-h) 
    & \scriptsize (NP-h) & \scriptsize (NP-h) &\scriptsize (NP-h) $^{**}$&\scriptsize (NP-h) $^{**}$
    \\ \hline

      \footnotesize split  & \scriptsize NP-h  &   \scriptsize NP-h &   \scriptsize NP-h    &     \scriptsize NP-h     &  \scriptsize NP-h $^{**}$ &  \scriptsize NP-h $^{**}$
  \\ 
   
 & \scriptsize (NP-h) & \scriptsize (NP-h)  & \scriptsize (NP-h)  &   \scriptsize (NP-h) & \scriptsize (NP-h) $^{**}$  & \scriptsize (NP-h) $^{**}$
 \\ \hline
    
   \footnotesize dually   & \scriptsize NP-h  & \scriptsize NP-h    & \scriptsize NP-h   &  \scriptsize NP-h & \scriptsize NP-h $^{**}$ & \scriptsize NP-h $^{**}$
   \\ 
  \footnotesize chordal  & \scriptsize (?) & \scriptsize (NP-h) $^*$ &  \scriptsize (NP-h) $^{**}$ & \scriptsize (NP-h) $^{**}$ & \scriptsize (NP-h) $^{**}$ & \scriptsize (NP-h) $^{**}$ 
  \\ \hline
  
  \footnotesize doubly    & \scriptsize NP-h & \scriptsize NP-h & \scriptsize NP-h 
&  \scriptsize NP-h & \scriptsize NP-h $^{**}$ & \scriptsize NP-h $^{**}$
  \\ 
  \footnotesize chordal  & \scriptsize (?) & \scriptsize (NP-h) $^{**}$  & \scriptsize (NP-h) $^{**}$ & \scriptsize (NP-h) $^{**}$ & \scriptsize (NP-h) $^{**}$& \scriptsize (NP-h) $^{**}$
  \\ \hline
   
  \footnotesize strongly  & \scriptsize P & \scriptsize P &     \scriptsize P & \scriptsize P & \scriptsize P $^{*}$ & \scriptsize P $^{**}$
  \\
  \footnotesize chordal  & \scriptsize (P) &  \scriptsize (P) & \scriptsize (P) & \scriptsize (P)  & \scriptsize (P) $^{*}$ & \scriptsize (P) \textbf{$^{**}$}
  \\ \hline
   
  \footnotesize bounded  & \scriptsize P & \scriptsize P & \scriptsize P  & \scriptsize ? & \scriptsize P $^{**}$ & \scriptsize ?
  \\ 
   \footnotesize clique-width & \scriptsize (P) & \scriptsize (P) $^{*}$& \scriptsize (P) $^{*}$ &\scriptsize (?)& \scriptsize (P) $^{**}$ & \scriptsize (?) 
   \\ \hline

\end{tabular}
\caption{
Doubly starred entries correspond to results established in this paper, while single stars denote new results derived from them. Question marks indicate unsolved problems
.}
\label{tablanueva}
\end{figure}

Some of the unknown computational complexity results are somewhat surprising, such as the computational complexity of $k$LP ($k \geq 2$) for bipartite, dually chordal, and doubly chordal graphs, especially if we consider that LP is known to be NP-hard in these families. On the other hand, we can say that $\mathcal L$DF, $\mathcal L$PF, GDF, and GPF are polynomial-time solvable for graphs with bounded clique-width only when restricted to instances with bounded function values while for general instances the computational complexity is unknown. 
Finally, although these problems were not addressed in this paper, the computational complexity of the Minus Domination problem remains unknown for dually and doubly chordal graphs. In contrast, the Signed Domination problem is known to be NP-hard in doubly chordal graphs. As noted in Section \ref{unificando}, these two problems exhibit significant similarities as specific instances of GDF. 
Indeed,
while Minus Domination problem is equivalent to GDF restricted to instances with $\u= 2 \uno$ and $k(v)= 2+ \delta(v)$, for all $v \in V$, Signed domination problem corresponds to instances of GDF with $\u=\uno$ and $k(v)=1+\lceil \frac{\delta(v)}{2}\rceil$.
On the \emph{packing side}, by Lemma \ref{iguales}, Signed Domination Problem is equivalent to GPF for instances with $\u=\uno$ and $k(v)= 
\lfloor \frac{\delta(v)}{2}\rfloor$ for all $v\in V$, and Minus Domination problem, to  GPF with $\u=2 \uno$ and $k(v)= \delta(v)$, for all $v\in V$.

In other line,
since the results on bounded clique-width graphs are based on Courcelle's Theorem (\cite{cour90}, \cite{coueng}), the algorithm that can be derived from this result does not take advantage of the structural characteristics of graphs belonging to some subfamilies.  
In \cite{lagos}, some technical results related to the behavior of the general packing number under various graph operations were studied allowing the design of specific algorithms to solve GPF on
block, spider, and
$P_4$-tidy graphs. 
This type of analysis can be extended, both for GPF and GDF, in order to obtain  specific algorithms and closed formulas for new families of graphs of bounded clique-width.

It is also interesting to note that strongly chordal graphs remain the only known graph class with unbounded clique-width for which the generalized packing and domination problems (GDF and GPF) are tractable. 

\nocite{*}

\bibliographystyle{itor}
\bibliography{itor}

\begin{thebibliography}{27}
\expandafter\ifx\csname natexlab\endcsname\relax\def\natexlab#1{#1}\fi
\expandafter\ifx\csname url\endcsname\relax
  \def\url#1{\texttt{#1}}\fi
\expandafter\ifx\csname urlprefix\endcsname\relax\def\urlprefix{URL }\fi
\providecommand{\eprint}[2][]{\url{#2}}
\bibitem[{Argiroffo et~al.(2015)Argiroffo, Leoni and Torres}]{ALT2015}
Argiroffo, G., Leoni, V., Torres, P., 2015.
\newblock On the complexity of {k}-domination and k-tuple domination in graphs.
\newblock \textit{Information Processing Letters} 115, 556--561.
\bibitem[{Argiroffo et~al.(2017)Argiroffo, Leoni and Torres}]{ALT}
Argiroffo, G., Leoni, V., Torres, P., 2017.
\newblock On the complexity of the labeled domination problem in graphs.
\newblock \textit{Int. Trans. in Op. Res.} 24, 355--367.
\bibitem[{Bertossi(1984)}]{Bertossi}
Bertossi, A., 1984.
\newblock Dominating sets for split and bipartite graphs.
\newblock \textit{Inform. Process. Lett.} 19, 37--40.
\bibitem[{Brandst\"adt et~al.(1998)Brandst\"adt, Chepoi and Dragan}]{bran}
Brandst\"adt, A., Chepoi, V.D., Dragan, F., 1998.
\newblock The algorithmic use of hypertree structure and maximum neighbourhood
  orderings.
\newblock \textit{Discrete Appl. Math.} 82, 43--77.
\bibitem[{Chang(1982)}]{Chang82}
Chang, G.J., 1982.
\newblock $k$-domination and graph covering problems.
\newblock Ph.D. thesis, School of OR and IE, Cornell University, Ithaca, NY.
\bibitem[{Chen(2001)}]{chen}
Chen, Y.M., 2001.
\newblock The fault tolerant domination problem on strongly chordal graphs.
\newblock Master's thesis, National Chung Cheng University, Taiwan.
\bibitem[{Courcelle(1990)}]{cour90}
Courcelle, B., 1990.
\newblock The monadic second-order logic of graphs. i. recognizable sets of
  nite graphs.
\newblock \textit{Information and Computation} 85, 1, 12--75.
\bibitem[{Courcelle et~al.(1993)Courcelle, Engelfriet and Rozenberg}]{coueng}
Courcelle, B., Engelfriet, J., Rozenberg, G., 1993.
\newblock Handle-rewriting hypergraph grammars.
\newblock \textit{J. Comput. System Sci.} 46, 218--270.
\bibitem[{Courcelle et~al.(2000)Courcelle, Makowsky and Rotics}]{cmr}
Courcelle, B., Makowsky, J.A., Rotics, U., 2000.
\newblock Linear time solvable optimization problems on graphs of bounded
  clique width.
\newblock \textit{Theory of Computing Systems} 33, 125--150.
\bibitem[{Dobson et~al.(2017)Dobson, Hinrichsen and Leoni}]{DHLITOR2016}
Dobson, M.P., Hinrichsen, E., Leoni, V., 2017.
\newblock On the complexity of the $\{k\}$-packing function problem.
\newblock \textit{International Transactions in Operational Research} 24,
  347--354.
\bibitem[{Dobson et~al.(2010)Dobson, Leoni and Nasini}]{ISCODLNpacking}
Dobson, M.P., Leoni, V., Nasini, G., 2010.
\newblock The $k$-limited packing and $k$-tuple domination problems in strongly
  chordal, $p_4$-tidy and split graphs.
\newblock \textit{Electronic Notes in Discrete Mathematics} 36, 559--566.
\bibitem[{Dobson et~al.(2011)Dobson, Leoni and Nasini}]{dln}
Dobson, M.P., Leoni, V., Nasini, G., 2011.
\newblock The multiple domination and limited packing problems in graphs.
\newblock \textit{Information Proccesing Letters} 111, 1108--1113.
\bibitem[{Farber(1983)}]{farber}
Farber, M., 1983.
\newblock Characterizations of strongly chordal graph.
\newblock \textit{Discrete Mathematics} 43, 173--189.
\bibitem[{Gallant et~al.(2010)Gallant, Gunther, Hartnell and Rall}]{GGHR1}
Gallant, R., Gunther, G., Hartnell, B., Rall, D., 2010.
\newblock Limited packing in graphs.
\newblock \textit{Discrete Applied Mathematics} 158, 1357--1364.
\bibitem[{Harary and Haynes(2000)}]{HH}
Harary, F., Haynes, T., 2000.
\newblock Double domination in graphs.
\newblock \textit{Ars Combin.} 55, 201--213.
\bibitem[{Haynes et~al.(1998)Haynes, Hedetniemi and Slater}]{Haynes}
Haynes, T., Hedetniemi, S., Slater, P., 1998.
\newblock \textit{Fundamentals of Domination in Graphs}.
\newblock Marcel Dekker, New York.
\bibitem[{He and Liang(2011)}]{kdom}
He, J., Liang, H., 2011.
\newblock Complexity of total {k}-domination and related problems.
\newblock \textit{Lecture Notes in Computer Science} 6681, 147--155.
\bibitem[{Hinrichsen and Leoni(2014)}]{HL}
Hinrichsen, E., Leoni, V., 2014.
\newblock $\{k\}$-packing functions in graphs.
\newblock \textit{Lecture Notes in Computer Science, Springer, Heidelberg}
  8596, 325--335.
\bibitem[{Hinrichsen et~al.(2019)Hinrichsen, Leoni and Safe}]{hls}
Hinrichsen, E., Leoni, V., Safe, M., 2019.
\newblock Labelled packing functions in graphs.
\newblock \textit{Electronic Notes in Information Processing Letters} 787.
\bibitem[{Hinrichsen et~al.(2023)Hinrichsen, Nasini and Vansteenkiste}]{lagos}
Hinrichsen, E., Nasini, G., Vansteenkiste, N., 2023.
\newblock On general packing functions in graphs.
\newblock \textit{Procedia Computer Science} 223.
\bibitem[{Lee and Chang(2008)}]{leechang}
Lee, C.M., Chang, M.S., 2008.
\newblock Variations of y-dominating functions on graphs.
\newblock \textit{Discrete Mathematics} 308, 4185--–4204.
\bibitem[{Leoni and Nasini(2014)}]{LNreductions}
Leoni, V., Nasini, G., 2014.
\newblock Limited packing and multiple domination problems: polynomial time
  reductions.
\newblock \textit{Discrete Applied Mathematics} 164, 547--553.
\bibitem[{Liao and Chang(2002)}]{liaochang}
Liao, C.S., Chang, G., 2002.
\newblock Algorithmic aspect of k-tuple domination in graphs.
\newblock \textit{Taiwanese Journal of Mathematics} 6, 3, 415--–420.
\bibitem[{Meir and Moon(1975)}]{MM}
Meir, A., Moon, J.W., 1975.
\newblock Relations between packing and covering numbers of a tree.
\newblock \textit{Pacific J. Math.} 61, 225--233.
\bibitem[{Müller and Brandstädt(1987)}]{MB}
Müller, H., Brandstädt, A., 1987.
\newblock The np-completeness of steiner tree and dominating set for chordal
  bipartite graphs.
\newblock \textit{Theoret. Comput. Sci.} 53, 257--265.
\bibitem[{Ore(1962)}]{ore}
Ore, O., 1962.
\newblock Theory of graphs.
\newblock \textit{American Mathematical Society Colloquium Publications} 38.
\bibitem[{Oum and Seymour(2006)}]{OumSey}
Oum, S., Seymour, P., 2006.
\newblock Approximating clique-width and branch-width.
\newblock \textit{Journal of Combinatorial Theory, Series B} 96, 514--528.

\end{thebibliography}

\section{Appendix} \label{apendice}

In this section, we provide the definitions for the problems mentioned in Section \ref{preliminary} and outline the results of computational complexity summarized in Table \ref{tabla}. Although most of these problems were initially formulated as decision problems, we provide their definitions in the optimization framework for consistency and clarity.
From now on, $G=(V,E)$ and $k\in \Z_+$.

\medskip

A dominating set of $G$ is a subset
$D$ of $V$ such that $|N[v]\cap D|\geq 1$ for all $v \in V$ and the dominating number $\gamma(G)$ is the minimum cardinality of a
dominating set of $G$. 
A $2$-packing of $G$ is a subset $B$ of $V$ where the distance between two vertices in $B$ is at least 3 \citep{MM}. This requirement is equivalent to $|N[v]\cap B|\leq 1$ for any $v\in V$, defined by \cite{GGHR1} as a $1$-limited packing of $G$.
The 2-packing number of $G$ is the
number of nodes in any largest 2-packing of $G$, denoted by $P_2(G)$. 
The related optimization problems are as follows: 

\medskip

{\bf Domination Problem} \citep{ore}

Instance: A graph $G$.

Objective: To obtain $\gamma(G)$.

\medskip

{\bf $2$-Packing Problem} \citep{MM}

Instance: A graph $G$.

Objective: To obtain $P_2(G)$.\\

These problems gave rise to the following variations:
a $k$-tuple dominating set (resp. $k$-limited packing) of $G$ is a subset of vertices $B\subset V$ that verify $|N[v]\cap B| \geq k$ (resp. $|N[v]\cap B| \leq k$) for every $v \in V$.
In addition, let $\gamma_{k}(G)$ (resp. $L_k(G)$) be the minimum (resp. maximum) cardinality of a $k$-tuple dominating set (resp. $k$-limited packing) of $G$.   

\medskip

For each fixed $k\geq 2$ the following problems are defined:

\medskip 

{\bf $k$LD: $k$-tuple Domination  Problem}  \citep{HH}

Instance: A graph $G$.

Objective: 
To obtain
$\gamma_{k}(G)$.
\medskip

{\bf $k$LP: $k$-limited Packing  Problem} \citep{GGHR1}

Instance: A graph $G$.

Objective: To obtain 
$L_{k}(G)$.

\medskip

The following two problems were subsequently defined, the above ones being particular instances of them:

\medskip 

{\bf LD: Multiple Domination Problem} \citep{dln}

Instance: A graph $G$, $k\geq 2$.

Objective: To obtain 
$\gamma_{k}(G)$.\\

\medskip 

{\bf LP: Limited Packing Problem} \citep{dln}

Instance: A graph $G$, $k\geq 2$.

Objective: To obtain 
$L_{k}(G)$.

\medskip 

Previous problems were generalized by allowing different quantities in the neighborhoods of each vertex and by requiring (or, respectively, forbidding) the presence of certain vertices in the selected sets.
Given $\k\in \Z_+^V$, a fault tolerant dominating set
is a subset $D$ of $V$ such that $|N[v]\cap D| \geq k(v)$ for
all $v \in V$.

\medskip 

{\bf Fault Tolerant Domination Problem} \citep{leechang}

Instance: A graph $G$, $\k \in \Z_+^V$.

Objective: To obtain 
a fault tolerant dominating set
of minimum cardinality of $G$.

\medskip 

Given a vector $\k  \in \Z^{V}_
+$ and $\mathcal R \subseteq V$ (resp. $\mathcal A \subseteq V$), a subset $D$ of vertices is a $(\k,\mathcal R)$-dominating set (resp. $({\k},\mathcal A)$-limited packing)  of $G$ if $R \subseteq D$ (resp. $D \subseteq \mathcal A$) and, for every $v \in V$,
$|N[v] \cap D| \geq k(v)$ (resp. $|N[v] \cap D| \leq k(v)$). From these definitions, the following problems are introduced:

\medskip 

{\bf GMD: Generalized Multiple Domination Problem} \citep{LNreductions}

Instance: A graph $G$, $\k \in \Z_+^V$, $\mathcal R \subseteq V$.  

Objective: To obtain the minimum size of a $(\k,\mathcal R)$-dominating set of $G$. 
\medskip

\medskip 

{\bf GLP: Generalized Limited Packing Problem} \citep{LNreductions}

Instance: A graph $G$, $\k \in \Z_+^V$, $\mathcal A \subseteq V(G)$.

Objective: To obtain the maximum size of a $({\k},\mathcal A)$-limited packing of $G$.

\medskip 

Subsequently, the following variations were defined.
A function $f : V \rightarrow [0, k]$ is a $\{k\}$-dominating function of $G$ (resp. $\{k\}$- packing function) if $f (N[v]) \geq k$ (resp. $f (N[v]) \leq k$) for all $v \in V$. Being $\gamma_{\{k\}}(G)$ (resp. $L_{\{k\}}(G)$) the minimum (resp. maximum) possible weight of a $\{k\}$- dominating (resp. $\{k\}$- packing) function of $G$, the following optimization problems are defined for each $k\geq 2$: 

\medskip 

{\bf $\{k\}$DF: $\{k\}$-dominating Function Problem} \citep{kdom}

Instance: A graph $G$.

Objective: To obtain 
$\gamma_{\{k\}}(G)$.

\medskip 

{\bf $\{k\}$PF: $\{k\}$-packing Function Problem } \citep{HL}

Instance: A graph $G$.

Objective: To obtain 
$L_{\{k\}}(G)$.

\medskip 

Given integers
$l, d, I$ with $l, d \geq 1$, let $Y = \{I, I+ d, I+ 2d, \ldots , I+ (l- 1)d\}$. A labelling function $L$ of $G$ is an assignment of a label $L(v) = (t(v), k(v))$ to each $v \in V$, where $t(v) \in \{F\} \cup Y$
and $k(v)$ is an integer. An $L$-dominating function of $G$ is a labelling function $f : V\rightarrow Y$ that satisfies the
following two conditions:
\begin{itemize}
    \item if $t(v) \neq F$, then $f (v) = t(v)$,
    \item $f (N[v]) \geq k(v)$ for all $v \in V$.
\end{itemize}
The weight of a $L$-dominating function $f$ is $f(V)$.

\medskip

{\bf{$\mathcal L$dom: Labelled Domination problem}} \citep{leechang} 

Instance: A graph $G$, a labelling function $ L$.

Objective: To obtain 
the minimum weight of an $L$-dominating function of $G$.

\medskip 

{\bf M-domination Problem (Mdom)} \citep{liaochang} is defined as $\mathcal L$dom restricted to instances where $Y = \{0, 1\}$ and $t(v) \in \{F, 1\}$. In addition, {\bf Signed}  and {\bf Minus Domination Problems} \citep{leechang} are defined as $\mathcal L$dom restricted to instances where $Y = \{-1, 1\}$ and $Y = \{-1,0,1\}$, respectively, and, in both cases,  $t(v)=F$ and $k(v)=1$, for all $v\in V$. 

\medskip 

Similarly,  labelled packing functions of $G$ are defined from labelling functions where $Y=[0,\ell]$ and $t(v)\in \{0,F\}$ for all $v\in V$. Given such a labelling function $L$, $f:V\rightarrow
Y$ is an $L$-packing function if 
$f(N[v])\leq k(v)$ and, if $t(v)\neq F$, then $f(v)=0$, for all $v\in V$.

\medskip 

{\bf{$\mathcal L$PF} : Labelled Packing problem} \citep{hls}

Instance: A graph $G$, a labelling function $L$ with $Y=[0,\ell]$ and $t(v)\in \{0,F\}$ for all $v\in V$.

Objective: To obtain 
the maximum weight of an $L$-packing function of $G$.

\medskip

Regarding the references for the computational complexity results used to construct Table \ref{tabla}, we mention the following:

The Domination problem is NP-hard on bipartite and split graphs \citep{MB,Bertossi} while on dually chordal graphs, it is polynomial-time solvable \citep{bran}.  
A maximum 2-packing of a dually chordal graph $G$ can be obtained in linear time \citep{bran}. 
For $k\geq 2$, $k$LD is NP-hard on split, bipartite \citep{liaochang}, and even on doubly chordal graphs \citep{leechang}.
Concerning packing problems, 
$k$LP is NP-hard on split graphs and, on bipartite graphs, LP is NP-hard \citep{dln}.
In addition, it is polynomial-time solvable on bounded clique-width graphs. \citep{HL}.

For $k\geq 2$, $\{k\}$DF is NP-hard on split and bipartite graphs \citep{liaochang}.  In addition, is polynomial-time solvable on dually chordal graphs \citep{leechang}.
In addition, 
$\{k\}$PF is polynomial-time solvable on dually chordal and bounded clique-width graphs \citep{DHLITOR2016,HL}.
and it is NP-hard on split and bipartite graphs \citep{DHLITOR2016}.

$\mathcal L$dom and $\mathcal L$PF are polynomial-time solvable on strongly chordal graphs \citep{leechang, hls}.  
In Corollary 15 in \cite{ALT}
it is proved that $\mathcal L$dom (for a fixed $\ell$) is  polynomial-time solvable on bounded clique-width graphs.

\end{document}